\documentclass[10pt,reqno]{amsart}
\usepackage{amsmath}
\usepackage{amssymb}
\usepackage{amsthm}

\newlength{\defbaselineskip}
\newcommand{\setlinespacing}[1]%
           {\setlength{\baselineskip}{#1 \defbaselineskip}}

\numberwithin{equation}{section}

\newtheorem{thm}{Theorem}[section]

\newtheorem{lem}[thm]{Lemma}
\newtheorem{prop}[thm]{Proposition}

\theoremstyle{definition}

\theoremstyle{remark}
\newtheorem{rem}[thm]{Remark}
\numberwithin{equation}{section}

\begin{document}
\title[Reconstruction of the initial data]
{Reconstruction of the initial data from the trace of the solutions on an infinite time cylinder of damped wave equations}

\author{Seongyeon Kim, Sunghwan Moon and Ihyeok Seo}

\thanks{This research was supported by RS-2023-00217116 (S. Moon) and NRF-2022R1A2C1011312 (I. Seo).}


\address{Department of Mathematics Education, Jeonju University, Jeonju 55069, Republic of Korea}
\email{sy\_kim@jj.ac.kr}

\address{Department of Mathematics, College of Natural Sciences, Kyungpook National University, Daegu 41566, Republic of Korea}
\email{sunghwan.moon@knu.ac.kr}

\address{Department of Mathematics, Sungkyunkwan University, Suwon 16419, Republic of Korea}
\email{ihseo@skku.edu}

\begin{abstract}
In this paper, we consider two types of damped wave equations: the weakly damped equation and the strongly damped equation. We recover the initial velocity from the trace of the solution on a space-time cylinder.
This inverse problem is related to Photoacoustic Tomography (PAT), a hybrid medical imaging technique. PAT is based on generating acoustic waves inside of an object of interest and one of the mathematical problem in PAT is reconstructing the initial velocity from the solution of the wave equation measured on the outside of object. Using the spherical harmonics and spectral theorem, we demonstrate a way to recover the initial velocity.
\end{abstract}

\maketitle

\section{Introduction}\label{sec1}
We consider the Cauchy problems for two distinct types of wave equations, which are the weakly damped equation \eqref{WW}, \eqref{WW1} and the strongly damped equation \eqref{SW}, \eqref{SW1}. The former is as follows:
\begin{align}\label{WW}
u_{tt} -\Delta u  +\gamma u_t &= 0
\,\,\,\qquad\text{on } \mathbb R^n \times [0,\infty),\\
\label{WW1}
u(x,0)=0, \quad
u_t(x,0) &= f(x) \quad\text{on }\mathbb R^n, 
\end{align}
where the term $\gamma u_t$ ($\gamma>0$) describes the weak damping effect,  acting to gradually reduce the energy of the wave. In the physical context, $\gamma$ can be seen as the friction coefficient that accounts for the damping effects of resistance.
The latter is as follows:
\begin{align}\label{SW}
u_{tt} -\Delta u  -\delta \Delta u_t &= 0\,\,\,\qquad\text{on }\mathbb R^n \times [0,\infty),\\
\label{SW1}
u(x,0)=0, \quad
u_t(x,0) &= f(x) \quad\text{on }\mathbb R^n,
\end{align}
where the term $\delta \Delta u_t$ ($\delta>0$) represents the strong damping effect, significantly reducing the energy of the wave. In this context, $\delta$ is a damping factor reflecting the combined influence of strain and strain rate in models of viscous material undergoing longitudinal vibrations \cite{CS,G,S,W}.

Let $\mathbb{S}^{n-1}$ denote the unit sphere in $\mathbb{R}^n$ and define the forward trace map 
\begin{equation}\label{trace}
\mathcal F: f\rightarrow u|_{\mathbb{S}^{n-1}\times [0,\infty)}.
\end{equation}
In this paper we study the inversion of $\mathcal{F}$.
More precisely, 
we are concerned with establishing a formula for reconstructing the initial velocity from the trace of the solutions on the infinite time cylinder $\mathbb{S}^{n-1}\times[0,\infty)$ to \eqref{WW}, \eqref{WW1} and \eqref{SW}, \eqref{SW1}. 

Thanks to the Fourier transform, solutions to linear PDEs can be generally formulated in terms of initial data.
It often emerges in practical applications that the reverse process,  recovering the unknown initial data from solutions measured at a detector surface, is of even greater importance. 
For example, one of the mathematical problems in Photoacoustic Tomography (PAT) boils down to recovering the initial data from the solution of the wave equation measured on a certain of a surface where the point-like detectors are located \cite{kuchmentk08,kuchment14book}. 

PAT, a hybrid imaging technique that combines high optical contrast with good ultrasound resolution, is based on generating ultrasound waves inside an object of interest by pulsed optical illumination.
The initial velocity of the induced sound waves contains biological information of the object, which is of great interest for medical diagnostics \cite{moonjmaa18,xuw05,zangerl19}. 
Therefore, the problem of reconstructing the initial velocity from data that are solutions to the wave equation measured at the detector surface outside the object is important.

In \cite{ESS,kowars10}, several wave equations considering attenuation and their causal properties are studied and the ill-posedness of the photoacoustic imaging problem in the case of an attenuating medium is studied.
\cite{FR,FPR} are devoted to the reconstruction of the initial data supported in a ball from the trace measured on a space-time cylinder.
In particular, \cite{FPR} studied this issue on a finite time cylinder, $\mathbb{S}^{n-1}\times [0,T]$ for a given $T>0$ and the odd dimension $n$; it's crucial to highlight that \cite{FPR} derived integral formulas for the inverse, contrasting with the series-based formulas presented in this paper, albeit still with constraints on the initial data's support (the $\gamma=0$, $\delta=0$ case). 
Also, in \cite{FR} the authors showed the uniqueness of the initial data from the solution of the wave equation on a finite time cylinder for all dimensions using the inversion formula for the spherical Radon transform.
Although we requires the trace of the solution on an infinite time cylinder, $\mathbb{S}^{n-1}\times [0,\infty)$ (because Huygens' principle does not work when the number of spatial dimensions is even), the problem considering damping effects ($\gamma> 0$ and $\delta> 0$) addressed in this paper has not been studied before, and our method differs from those in previously relevant works.

Before stating our results, we introduce some notations.
For $g\in L^2(\mathbb{R}^n)$, $g_{lk}$ is the coefficient of the spherical harmonics expansion (see \cite{S1,SW})
$$g(y)= \sum_{l=0}^{\infty} \sum_{k=0}^{d(n,l)}  g_{lk}(\rho){\bf{Y}}_{lk} (\omega), \quad y=\rho \omega, \quad\omega\in \mathbb{S}^{n-1},\quad\rho>0,$$
where the spherical harmonics ${\bf Y}_{lk}$ form an orthonormal basis in $L^2(\mathbb{S}^{n-1})$, $d(n,l)=(n+l-3)!(2l+n-2)/(l!(n-2)!)$ for $l \in \mathbb{N}$, and $d(n,0)=1$.
We also denote by $m(t,|\xi|)$ the Fourier multiplier appeared in the solution formula 
\begin{equation}\label{s}
u(x,t)=	\int_{\mathbb{R}^n} e^{ix\cdot \xi} m(t,|\xi|) \hat f(\xi) d\xi,
\end{equation}
where the exact form of $m$ will be detailed further in Section \ref{sec2} (refer to \eqref{mull} for \eqref{SW}, \eqref{SW1}, and \eqref{mull2} for \eqref{WW}, \eqref{WW1}).
Finally, $\bf{H}$ is the integral operator on $L^2([0,\infty)])$ given by 
\begin{equation*}
{\bf{H}}(h)(s)=\int_0^{\infty} H(s,\rho)h(\rho)d\rho
\end{equation*} 
where $H(s,\rho)=\int_0^{\infty}K(t,s)K(t,\rho)dt$ with $K(t,\rho)=m(t,\rho)\rho^{\alpha}J_{l+\frac{n-2}2}(\rho).$
Here, $J_{\nu}$ denotes the Bessel function of the order $\nu$ with $\text{Re}\nu>-1/2.$ 

We will prove that ${\bf H}$ is a compact self-adjoint operator on the infinite-dimensional separable Hilbert space $L^2([0,\infty))$. See below Lemma \ref{ker}.
Then by Theorem \ref{Sthm} (and Remark \ref{rem}), there exists an orthonormal basis $\{v_h\}$ 
for $(\text{ker}\,{\bf{H}})^{\perp}$ consisting of eigenfunctions of ${\bf{H}}$, with corresponding nonzero eigenvalues $\lambda_h$.

Let $f_{lk}$ and $u_{lk}$ be the coefficients of the spherical harmonics expansions of $f$ and $u$, respectively. Now we state our results which provide a formula for reconstructing the initial velocity from the trace of the solution on the infinite time cylinder $\mathbb{S}^{n-1}\times[0,\infty)$:

\begin{thm}\label{thm}
Let $n\ge2$. Let $u\in C([0,\infty);L^2(\mathbb R^n))$ be the solution to \eqref{WW}, \eqref{WW1} with $f\in L^2(\mathbb{R}^n)$, restricted on $\mathbb{S}^{n-1}$.
Then we have 
\begin{equation}\label{f}
\hat f_{lk}(\rho) = (2\pi)^{-\frac{n}{2}}i^{l}\rho^{\alpha-\frac{n}{2}}\sum_{h=1}^\infty \frac{1}{\lambda_h}v_h(\rho)\int_0^\infty v_h(s)\int_0^{\infty}m(t,s)s^\alpha J_{l+\frac{n-2}{2}}(s)u_{lk}(t)dtds 
\end{equation}
if $\hat f_{lk}(\rho) \rho^{\frac{n}{2}-\alpha}\in L^2([0,\infty))\cap(\text{ker}\,{\bf{H}})^{\perp}$ for some $-\frac{n-3}{2}<\alpha<1$.
\end{thm}

\begin{thm}\label{thm2}
Let $n\ge2$. Let $u\in C([0,\infty);L^2(\mathbb R^n))$ be the solution to \eqref{SW}, \eqref{SW1} with $f\in L^2(\mathbb{R}^n)$, restricted on $\mathbb{S}^{n-1}$. 
Then we have the same reconstruction formula \eqref{f} 
if $\hat f_{lk}(\rho) \rho^{\frac{n}{2}-\alpha}\in L^2([0,\infty))\cap (\text{ker}\,{\bf{H}})^{\perp}$ for some $-\frac{n-5}{2}<\alpha<2$.
\end{thm}

\begin{rem}
The restriction of $u$ on $\mathbb{S}^{n-1}\times[0,\infty)$ in the theorems is plausible in view of the fact that the point-like detectors are located on the unit sphere.
\end{rem}

Finally, we note that an explicit condition on $f$ for which 
$\hat f_{lk}(\rho) \rho^{\frac{n}{2}-\alpha}\in L^2([0,\infty))$
can be given as $|\nabla|^{(1-2\alpha)/2}f\in L^2(\mathbb{R}^n)$.
Indeed, if $|\nabla|^{(1-2\alpha)/2} f\in L^2(\mathbb{R}^n)$, then
$|\xi|^{\frac{1-2\alpha}2}\hat f (\xi)\in L^2(\mathbb{R}^n)$ by Plancherel's theorem.
Now, we write $\hat{f}$ using the spherical harmonics expansion:
$$\hat f(\rho\omega)=\sum_{l=0}^\infty\sum_{k=0}^{d(n,l)} \hat f_{lk}(\rho) {\bf{Y}}_{lk}(\omega).$$
Then, 
\begin{align*}
\big\||\xi|^{\frac{1-2\alpha}2}\hat f (\xi)\big\|_{L^2(\mathbb{R}^n)}
&=\int_0^\infty\int_{\mathbb S^{n-1}}|\rho^{\frac{1-2\alpha}{2}}\hat f(\rho\omega)|^2d\omega \rho^{n-1}d\rho\\
&=\int_0^\infty\int_{\mathbb S^{n-1}}\Big|\sum_{l=0}^\infty \sum_{k=0}^{d(n,l)} \rho^{\frac{1-2\alpha}{2}}\hat f_{lk}(\rho) {\bf{Y}}_{lk}(\omega)\Big|^2 d\omega \rho^{n-1}d\rho\\
&=\int_0^\infty \sum_{l=0}^\infty\sum_{k=0}^{d(n,l)} \Big|\rho^{\frac{1-2\alpha}{2}}\hat f_{lk}(\rho)\Big|^2 \rho^{n-1}d\rho
\end{align*}
by the orthonormality of ${\bf{Y}}_{lk}$ on $L^2(S^{n-1})$.
Consequently, it follows that $\hat f_{lk}(\rho)\rho^{n/2-\alpha}\in L^2([0,\infty))$.
Note also that if $\alpha$ is taken as $1/2$ (which is possible when $n\geq 3$ for Theorem \ref{thm} and $n\geq 5$ for Theorem \ref{thm2}), the condition $|\nabla|^{(1-2\alpha)/2}f \in L^2(\mathbb{R}^n)$ changes to $f \in L^2(\mathbb{R}^n)$, no longer imposing any additional conditions.

\

\noindent\textit{Outline of the paper.}
In Section \ref{sec2}, we represent the solutions by \eqref{s} utilizing the Fourier transform. To obtain the reconstruction formula \eqref{f} in Theorems \ref{thm} and \ref{thm2}, we substitute the spherical harmonics expansion of $\hat f$ into \eqref{s} and compare the resulting coefficients with $u_{lk}$. This leads us to 
\begin{equation*}
u_{lk}(t)=(2\pi)^{\frac{n}{2}}i^{-l}\int_0^{\infty} K(t,\rho) \hat f_{lk}(\rho)\rho^{\frac{n}{2}-\alpha} d\rho.
\end{equation*}
Multiplying both sides of this by $K(t,s)$ and integrating with respect to $t$, we derive the following integral equation of the first kind
\begin{equation}\label{H1}
\int_0^{\infty} K(t,s) u_{lk}(t)dt=\int_0^{\infty} H(s,\rho) \big((2\pi)^{\frac{n}{2}}i^{-l}\hat f_{lk}(\rho)\rho^{\frac{n}{2}-\alpha}\big) d\rho.
\end{equation}
From the compactness and self-adjointness of $\bf{H}$, we use an orthonormal basis $\{v_h\}$ for $(\text{ker}\,{\bf{H}})^{\perp}$ consisting of eigenfunctions of $\bf H$ in $L^2([0,\infty))$ to compare the coefficients in the eigenfunction expansion of each side of \eqref{H1}. 
This leads us to the reconstruction formula \eqref{f}. The above is the main content of Section \ref{sec3}.
In the final section, Section \ref{sec4}, we prove $H(s,\rho)\in L^2([0,\infty)) \times L^2([0,\infty))$ to guarantee that $\bf{H}$ is a compact self-adjoint operator on $L^2([0,\infty))$.

\section{Representation of the solution}\label{sec2}
In this section, we present the solutions $u\in C([0,\infty);L^2(\mathbb{R}^n))$ to the damped wave equations \eqref{WW}, \eqref{WW1} and \eqref{SW}, \eqref{SW1} with $f\in L^2(\mathbb{R}^n)$ in terms of the Fourier transform.
By taking the Fourier transform in the spatial variable, we convert \eqref{SW}, \eqref{SW1} into the following ODE
\begin{equation}\label{ode}
	\begin{cases}
\partial_t^2 \hat u (\xi, t)+ |\xi|^2 \hat u (\xi, t)+\delta|\xi|^2 \partial_t \hat u(\xi, t) = 0,\\
\hat u (\xi,0)=0,\\
\partial_t \hat u(\xi,0)=\hat f(\xi),
\end{cases}
\end{equation}
in which roots of the characteristic equation
$$\mu^2 + \delta |\xi|^2 \mu +|\xi|^2=0$$
are complex roots in the form 
\begin{equation*} 
\mu=\begin{cases}
\frac{-\delta|\xi|^2 \pm |\xi|\sqrt{\delta^2|\xi|^2-4}}{2} \quad \textrm{if}\quad|\xi|>2/\delta,\\
\frac{-\delta|\xi|^2 \pm i|\xi|\sqrt{4-\delta^2|\xi|^2}}{2} \quad \textrm{if}\quad|\xi|<2/\delta.
\end{cases}
\end{equation*}
Here, we denote the roots as $\mu_1$ and $\mu_2$, without significance to their order.
Applying the initial condition, the solution to the ODE \eqref{ode} is given by 
\begin{align*}
&\hat u (\xi, t) = \frac{1}{\mu_1-\mu_2}\big(e^{\mu_1 t}-e^{\mu_2 t}\big)\hat f(\xi)\\
&\,\,=\begin{cases} \frac1{|\xi|\sqrt{\delta^2|\xi|^2 - 4}} \big(e^{\frac{-\delta|\xi|^2+|\xi|\sqrt{\delta^2|\xi|^2-4}}{2}t}-e^{\frac{-\delta|\xi|^2-|\xi|\sqrt{\delta^2|\xi|^2-4}}{2}t}\big)\hat f(\xi) \quad \textnormal{if} \quad |\xi|>2/\delta, \\
\frac2{|\xi|\sqrt{4-\delta^2|\xi|^2}}e^{-\frac{\delta|\xi|^2 t}{2}} \sin\big( \frac{t|\xi|\sqrt{4-\delta^2|\xi|^2}}{2}\big)\,\hat f(\xi) \quad \textnormal{if} \quad |\xi|<2/\delta.
\end{cases}
\end{align*}
Now taking the inverse Fourier transform, the solution to \eqref{SW}, \eqref{SW1} is expressed as
\begin{equation}\label{sol}
u(x,t)=\int_{\mathbb{R}^n}e^{ix\cdot\xi} \,m(t,|\xi|)\hat f(\xi)d\xi,
\end{equation}
where $m$ is the bounded Fourier multiplier given by
\begin{align}\label{mull}
m(t,|\xi|)&=\frac{1}{\mu_1-\mu_2}\big(e^{\mu_1 t}-e^{\mu_2 t}\big)\\
&=\begin{cases} \frac1{|\xi|\sqrt{\delta^2|\xi|^2 - 4}} \big(e^{\frac{-\delta|\xi|^2+|\xi|\sqrt{\delta^2|\xi|^2-4}}{2}t}-e^{\frac{-\delta|\xi|^2-|\xi|\sqrt{\delta^2|\xi|^2-4}}{2}t}\big) \quad \textnormal{if} \quad |\xi|>2/\delta, \\
\nonumber
\frac2{|\xi|\sqrt{4-\delta^2|\xi|^2}}e^{-\frac{\delta|\xi|^2 t}{2}} \sin \big(\frac{t|\xi|\sqrt{4-\delta^2|\xi|^2}}{2}\big) \quad \textnormal{if} \quad |\xi|<2/\delta.
	\end{cases}
\end{align}

In the same way, the solution to the weakly damped equation \eqref{WW}, \eqref{WW1} is given by \eqref{sol} with 
\begin{equation}\label{mull2}
m(t,|\xi|)=\begin{cases} \frac1{\sqrt{\gamma^2-4|\xi|^2}} \big(e^{\frac{-\gamma+\sqrt{\gamma^2-4|\xi|^2}}{2}t}-e^{\frac{-\gamma-\sqrt{\gamma^2-4|\xi|^2}}{2}t}\big) \quad \textnormal{if} \quad |\xi|<\gamma/2, \\
		\frac{2}{\sqrt{4|\xi|^2-\gamma^2}}e^{-\frac{\gamma t}{2}} \sin \big(\frac{t\sqrt{4|\xi|^2-\gamma^2}}{2}\big) \quad \textnormal{if} \quad |\xi|>\gamma/2. 
	\end{cases}
\end{equation}

Finally, the continuity of the solutions follows from applying the mean value theorem to $m(t,|\xi|)$ as
\begin{align*}
\|u(x,t)-u(x,t_0)\|_{L_x^2} &= \|\hat u(\xi,t)-\hat u(\xi,t_0)\|_{L_\xi^2} \\
&\leq \big\|m(t,|\xi|)-m(t_0,|\xi|)\big\|_{L_\xi^\infty}\|\hat f\|_{L_\xi^2} \\
&\leq C|t-t_0| \rightarrow 0 \quad \textrm{as} \,\,\, t\rightarrow t_0.
\end{align*}

\section{Reconstruction of the initial data}\label{sec3}
Now we prove Theorem \ref{thm2}, which shows the process of reconstructing the unknown initial function $f$ from the trace of the solution $u$ on the infinite time cylinder $\mathbb{S}^{n-1} \times [0, \infty)$
of the strongly damped equation \eqref{SW}, \eqref{SW1}. 
Since the multiplier $m(t,|\xi|)$ for the strongly damped case is more complex, we only provide the proof for this case.
(The other case is proved in the same way.)

By the spherical harmonics expansion of the function $\hat f$, we see that
$$\sum_{l=0}^{N}\sum_{k=0}^{d(n,l)} {\hat f}_{lk}(\rho) {\bf{Y}}_{lk}(\omega)\rightarrow\hat f(\xi)$$ 
in $L^2$ as $N\rightarrow\infty$. Since $m(t,|\xi|)$ is bounded on $\xi$, the Plancherel theorem now gives 
$$\int_{\mathbb{R}^n}e^{ix\cdot\xi} \,m(t,|\xi|)\sum_{l=0}^{N}\sum_{k=0}^{d(n,l)} {\hat f}_{lk}(\rho) {\bf{Y}}_{lk}(\omega)d\xi\rightarrow\int_{\mathbb{R}^n}e^{ix\cdot\xi} \,m(t,|\xi|)\hat f(\xi)d\xi$$ 
in $L^2$ as $N\rightarrow\infty$.
Using \eqref{sol} and polar coordinates $\xi=\rho \omega$ with $\omega \in \mathbb{S}^{n-1}$, $\rho>0$, 
we can rewrite the solution $u$ as 
\begin{align}\label{u}
\nonumber
u(\theta,t)&=\sum_{l=0}^{\infty}\sum_{k=0}^{d(n,l)}\int_0^{\infty} m(t,\rho)\hat f_{lk}(\rho) \rho^{n-1} \int_{\mathbb{S}^{n-1}} e^{i\theta\cdot \rho \omega} {\bf{Y}}_{lk}(\omega) d\omega d\rho\\
&=\sum_{l=0}^{\infty}\sum_{k=0}^{d(n,l)}(2\pi)^{\frac{n}2}i^{-l} \int_0^{\infty} m(t,\rho) \hat f_{lk}(\rho)\rho^{\frac{n}2} J_{l+\frac{n-2}{2}}(\rho) d\rho {\bf{Y}}_{lk}(\theta),\quad \theta\in\mathbb{S}^{n-1}.
\end{align}
Recall here that we are considering $u$ restricted on the unit sphere $\mathbb{S}^{n-1}$, and for the last equality we used the following identity
$$\int_{\mathbb{S}^{n-1}}e^{i\xi\cdot \omega} {{\bf{Y}}_{lk}(\omega)d\omega}=(2\pi)^{\frac{n}2}i^{-l}|\xi|^{-\frac{n-2}2}J_{l+\frac{n-2}{2}}(|\xi|){\bf{Y}}_{lk}\big(\frac{\xi}{|\xi|}\big)$$
where $J_\nu$ denotes the Bessel function of complex order $\nu$ with $\textrm{Re} \nu >-1/2$ (see \cite[Theorem 3.10, p. 158]{SW} and \cite[(2.9)]{WF}).
We then expand $u$ in terms of spherical harmonics as 
$$u(\theta,t)=\sum_{l=0}^{\infty}\sum_{k=0}^{d(n,l)} u_{lk}(t){\bf{Y}}_{lk}(\theta).$$
Since ${\bf{Y}}_{lk}$ form an orthonormal basis in $L^2(\mathbb{S}^{n-1})$, comparing this expansion with \eqref{u} gives 
\begin{equation}\label{b}
u_{lk}(t)=(2\pi)^{\frac{n}2}i^{-l}	\int_0^{\infty} m(t,\rho) \hat f_{lk}(\rho)\rho^{\frac{n}2} J_{l+\frac{n-2}{2}}(\rho) d\rho,
\end{equation}
where
\begin{equation*}
	m(t,\rho)=\begin{cases}
		\frac1{\rho\sqrt{\delta^2\rho^2 - 4}} \big(e^{\frac{-\delta\rho^2+\rho\sqrt{\delta^2\rho^2-4}}{2}t}-e^{\frac{-\delta\rho^2-\rho\sqrt{\delta^2\rho^2-4}}{2}t}\big) \quad \textnormal{if} \quad \rho>2/\delta, \\
		\frac2{\rho\sqrt{4-\delta^2\rho^2}}e^{\frac{-\delta\rho^2t}{2}} \sin \big(\frac{t\rho\sqrt{4-\delta^2\rho^2}}{2}\big) \quad \textnormal{if} \quad \rho<2/\delta.
	\end{cases}
\end{equation*}

We now note that \eqref{b} is the Fredholm integral equation of $(2\pi)^{\frac{n}2}i^{-l}\hat f_{lk}(\rho)\rho^{\frac{n}{2}-\alpha}$ with the kernel 
\begin{equation}\label{K1}
K(t,\rho)= m(t,\rho) \rho^{\alpha} J_{l+\frac{n-2}{2}}(\rho).
\end{equation}
Multiplying both sides of \eqref{b} by $K(t,s)$ and integrating with respect to $t$, we have
\begin{align}\label{H}
\nonumber
\int_0^{\infty}K(t,s) u_{lk}(t) dt&=(2\pi)^{\frac{n}2}i^{-l}\int_0^{\infty}\int_0^\infty K(t,s)K(t,\rho)(\hat f_{lk}(\rho)\rho^{\frac{n}{2}-\alpha})d\rho dt\\
&=(2\pi)^{\frac{n}2}i^{-l}\int_0^{\infty}\Big(\int_0^\infty K(t,s)K(t,\rho)dt\Big) (\hat f_{lk}(\rho)\rho^{\frac{n}{2}-\alpha})d\rho.
\end{align}
Letting 
\begin{equation}\label{g}
\int_0^{\infty}K(t,s)u_{lk}(t)dt=g_{lk}(s) \quad \textrm{and}\quad \int_0^{\infty} K(t,s)K(t,\rho)dt=H(s,\rho),
\end{equation}
the equation \eqref{H} is transformed into   
\begin{equation}\label{F}
g_{lk}(s)=(2\pi)^{\frac{n}2}i^{-l}\int_0^{\infty}H(s,\rho)(\hat f_{lk}(\rho)\rho^{\frac{n}{2}-\alpha})d\rho={\bf{H}}\big((2\pi)^{\frac{n}2}i^{-l}\hat f_{lk}(\rho)\rho^{\frac{n}{2}-\alpha}\big)(s).
\end{equation}

To make use of the following theorem, Theorem \ref{Sthm} (see, for example, \cite[Corollarly 5.4]{C}), with $T={\bf{H}}$, we shall prove that ${\bf H}$ is a compact self-adjoint operator on the infinite-dimensional separable Hilbert space $L^2([0,\infty))$.
We will postpone the proof of this for now, and continue with the proof of Theorem \ref{thm2}.
\begin{thm}\label{Sthm}
Let $T$ be a compact, self-adjoint operator on an infinite-dimensional separable Hilbert space $H$. Then there is a sequence $\{\lambda_h\}$ of real numbers and an orthonormal basis $\{v_h\}$ for $(\text{ker}\,T)^{\perp}$ such that for all $x\in H$,
$$Tx=\sum_{h=1}^\infty \lambda_h \big\langle x,v_h \big\rangle_H\, v_h,$$
where $\langle x,v_h \rangle_{H}$ denotes the inner product in the space $H$.
\end{thm}

\begin{rem}\label{rem}
Since $\{v_h\}$ is an orthonormal basis for $(\text{ker}\,T)^{\perp}$, every $\lambda_h$ is nonzero eigenvalue from the fact that
$$0\neq Tv_m=\sum_{h=1}^\infty \lambda_h \big\langle v_m , v_h\big\rangle_H\, v_h=\lambda_m v_m.$$
\end{rem}

By Theorem \ref{Sthm}, there exists an orthonormal basis $\{v_h\}$ 
for $(\text{ker}\,{\bf{H}})^{\perp}$ consisting of eigenfunctions of ${\bf{H}}$, with corresponding nonzero eigenvalues $\lambda_h$, and
\begin{equation}\label{zeroeig}
{\bf{H}} ((2\pi)^{\frac{n}2}i^{-l}\hat f_{lk}(\rho)\rho^{\frac{n}{2}-\alpha})(s)
=\sum_{h=1}^{\infty}\lambda_h\big\langle (2\pi)^{\frac{n}2}i^{-l}\hat f_{lk}(\rho)\rho^{\frac{n}{2}-\alpha}, v_h(\rho)\big\rangle_{L^2([0,\infty))}v_h(s).
\end{equation}
On the other hand,  
\begin{equation}\label{notk}
g_{lk}(s)=\sum_{h=1}^{\infty} \big\langle g_{lk} , v_h \big\rangle_{L^2([0,\infty))} v_h(s)
\end{equation}
since $g_{lk}\in (\text{ker}\,{\bf{H}})^{\perp}$ for each $l,k$.
Indeed, for $h\in\text{ker}\,{\bf{H}}$, note that
\begin{align*}
\big\langle g_{lk},h\big\rangle_{L^2([0,\infty))}&=\big\langle {\bf{H}}\big((2\pi)^{\frac{n}2}i^{-l}\hat f_{lk}(\rho)\rho^{\frac{n}{2}-\alpha}\big),h\big\rangle_{L^2([0,\infty))}\\
&=\big\langle (2\pi)^{\frac{n}2}i^{-l}\hat f_{lk}(\rho)\rho^{\frac{n}{2}-\alpha},{\bf{H}} h\big\rangle_{L^2([0,\infty))}
=0.
\end{align*}

From \eqref{zeroeig}, \eqref{notk} and \eqref{F}, we now deduce 
\begin{equation}\label{hj}
\big\langle g_{lk}, v_h \big\rangle_{L^2([0,\infty))} = \lambda_h \big\langle (2\pi)^{\frac{n}2}i^{-l}\hat f_{lk}(\rho)\rho^{\frac{n}{2}-\alpha}, v_h(\rho)\big\rangle_{L^2([0,\infty))}.
\end{equation}
Since $\hat f_{lk}(\rho) \rho^{\frac{n}2-\alpha}\in (\text{ker}\,{\bf{H}})^{\perp}$, as assumed in the theorem, we can write
\begin{align*}
(2\pi)^{\frac{n}2}i^{-l}\hat f_{lk}(\rho) \rho^{\frac{n}2-\alpha}&= \sum_{h=1}^{\infty} \big\langle (2\pi)^{\frac{n}2}i^{-l}\hat f_{lk}(\tilde\rho) \tilde\rho^{\frac{n}2-\alpha}, v_h(\tilde\rho) \big\rangle_{L^2([0,\infty))} v_h(\rho)\\
& = \sum_{h=1}^{\infty} \frac{1}{\lambda_h} \big\langle g_{lk}, v_h \big\rangle_{L^2([0,\infty))} v_h(\rho)\\ 
&= \sum_{h=1}^\infty \frac{1}{\lambda_h} \Big\langle \int_0^{\infty}K(t,s)u_{lk}(t)dt, v_h(s) \Big\rangle_{L^2([0,\infty))} v_h(\rho),
\end{align*}
where we used \eqref{hj} and the first in \eqref{g} for the second and last equalities, respectively. 
Substituting \eqref{K1} into the right-hand side immediately implies \eqref{f} as desired.

It remains to prove that ${\bf H}$ is a compact self-adjoint operator on $L^2([0,\infty))$.
For this we first need to obtain the following lemma which will be obtained in the next section.
\begin{lem}\label{ker}
For $n\ge 2$ and $-(n-5)/2<\alpha<2$, we have
\begin{equation}\label{kerr}
\int_0^\infty \int_0^\infty \big|K(t,\rho)\big|^2 dtd\rho<\infty,
\end{equation}
where $K(t,\rho)$ is given as in \eqref{K1}.
\end{lem}
\begin{rem}
For the weakly damped case, \eqref{kerr} holds for $n\ge2$ and $-(n-3)/2<\alpha<1$.
\end{rem}
Then, using the Cauchy-Schwarz inequality and the lemma, we have
\begin{align*}
\|H\|_{L_{s,\rho}^2}^2&=\Big\|\int_0^{\infty} K(t,s)K(t,\rho) dt\Big\|_{L_{s,\rho}^2}^2 \\
&\leq \big\|\|K(t,s)\|_{L_t^2}\|K(t,\rho)\|_{L_t^2}\big\|_{L_{s,\rho}^2}^2\\
&=\|K\|_{L_{t,\rho}^2}^4 <\infty,
\end{align*}
which implies that $\bf H$ is the compact operator by applying Proposition \ref{C} with ${\bf K} = {\bf H}$ and $k=H$.
\begin{prop}\cite[Proposition 4.7, p. 43]{C}\label{C}
	If $k\in L^2([0,\infty)\times[0,\infty))$, then 
	$${\bf{K}}f(x)=\int_0^{\infty} k(x,y)f(y)dy$$
	is a bounded compact operator on $L^2([0,\infty))$.
\end{prop}
The self-adjointness is trivial; 
if we set $H^\ast(s,\rho)=\overline{H(\rho,s)}$, then ${\bf{H^\ast}}$ is the adjoint operator with the kernel $H^\ast$. 
Since $H(\rho,s)=H(s,\rho)$ and the kernels are real-valued, $H^\ast(s,\rho)=H(s,\rho)$.
That is, $\bf{H}$ is self-adjoint.

\section{Proof of Lemma \ref{ker}}\label{sec4}
This section is devoted to proving Lemma \ref{ker}.
We first recall that 
\begin{equation*}
	K(t,\rho)= m(t,\rho)\rho^\alpha J_{l+\frac{n-2}{2}}(\rho), 
\end{equation*}
where \begin{equation}\label{mul}
	m(t,\rho)=\begin{cases}
		\frac1{\rho\sqrt{\delta^2\rho^2 - 4}} \big(e^{\frac{-\delta\rho^2+\rho\sqrt{\delta^2\rho^2-4}}{2}t}-e^{\frac{-\delta\rho^2-\rho\sqrt{\delta^2\rho^2-4}}{2}t}\big) \quad \textnormal{if} \quad \rho>2/\delta, \\
		\frac2{\rho\sqrt{4-\delta^2\rho^2}}e^{\frac{-\delta\rho^2t}{2}} \sin \big(\frac{t\rho\sqrt{4-\delta^2\rho^2}}{2} \big)\quad \textnormal{if} \quad \rho<2/\delta.
	\end{cases}
\end{equation}
To bound the integral
\begin{equation}\label{K}
	\int_0^{\infty}\int_0^{\infty} |m(t,\rho)J_{l+\frac{n-2}{2}}(\rho)|^2 \rho^{2\alpha} dtd\rho,
\end{equation}
we divide cases into two parts, $\delta>2$ and $\delta\leq2$.

\subsection{The case $\delta>2$}
Motivated the following fact (see, for example, \cite[Appendix B]{Ga}) that for $\textrm{Re} \nu>-1/2$
\begin{equation}\label{Be}
	|J_\nu (m)|\leq 	\begin{cases}
		C_\nu m^{\textrm{Re} \nu} \quad \textrm{if}\quad  0<m<1,\\
		C_\nu m^{-\frac12} \quad \textrm{if}\quad m\ge1,
	\end{cases}
\end{equation}
we first split the integral with respect to $\rho$ in \eqref{K} into the regions $0<\rho<2/\delta$, $2/\delta<\rho<1$ and $\rho>1$ as
\begin{align}\label{in}
&\int_0^{2/\delta}\int_0^{\infty} |m(t,\rho)|^2 \rho^{2\alpha}|J_{l+\frac{n-2}{2}}(\rho)|^2dtd\rho\\
\nonumber
&+\int_{2/\delta}^{1}\int_0^{\infty}  |m(t,\rho)|^2 \rho^{2\alpha}|J_{l+\frac{n-2}{2}}(\rho)|^2dt d\rho
+\int_{1}^{\infty}\int_0^{\infty}  |m(t,\rho)|^2 \rho^{2\alpha}|J_{l+\frac{n-2}{2}}(\rho)|^2dt d\rho.
\end{align}

Using \eqref{mul} and \eqref{Be}, the first part is then bounded by
\begin{equation}\label{f1}
\int_0^{2/\delta}\frac{4\rho^{2\alpha+2\l+n-2}}{\rho^2(4 -\delta^2 \rho^2)}\int_0^{\infty} e^{-\delta \rho^2 t} \sin^2 \big(\frac{t\rho\sqrt{4-\delta^2 \rho^2}}{2}\big)dt d\rho.
\end{equation}
Applying the integration by parts in $t$ after using $\sin^2x = \frac{1-\cos2x}{2}$, we have 
\begin{align}\label{int}
\nonumber
\int_0^{\infty} e^{-\delta\rho^2t} \sin^2\big({\frac{t\rho\sqrt{4-\delta^2\rho^2}}{2}}\big)dt
&=\frac{1}{2}\big(\int_0^{\infty} e^{-\delta\rho^2 t} dt-\int_0^{\infty} e^{-\delta\rho^2 t} \cos (t\rho\sqrt{4-\delta^2\rho^2})dt\big)\\
&=\frac12\big(\frac1{\delta\rho^2}-A\big),
\end{align}
where $A$ denotes the last integral.
We apply the integration by parts twice to calculate
\begin{align}
\nonumber
&A=\int_0^{\infty}e^{-\delta\rho^2 t} \cos \big(t\rho\sqrt{4-\delta^2\rho^2}\big)dt\\
\nonumber
&=\Big[-\frac1{\delta\rho^2}e^{-\delta\rho^2t}\cos \big(t\rho\sqrt{4-\delta^2\rho^2}\big)\Big]_0^{\infty} -\int_0^{\infty}
\frac{\rho\sqrt{4-\delta^2\rho^2}}{\delta\rho^2}e^{-\delta\rho^2t}\sin\big( t\rho\sqrt{4-\delta^2\rho^2}\big)dt\\
\nonumber
&=\frac1{\delta\rho^2}-\frac{\rho\sqrt{4-\delta^2\rho^2}}{\delta\rho^2}\int_0^{\infty} e^{-\delta\rho^2t}\sin\big(t\rho\sqrt{4-\delta^2\rho^2}\big)dt\\
\nonumber
&=\frac1{\delta\rho^2}-\frac{\rho\sqrt{4-\delta^2\rho^2}}{\delta\rho^2}\Big[-\frac1{\delta\rho^2}e^{-\delta\rho^2t}\sin\big(t\rho\sqrt{4-\delta^2\rho^2}\big)\Big]_0^{\infty}\\
\nonumber
&\qquad \qquad \qquad \qquad \qquad \qquad \qquad \quad -\frac{\rho^2(4-\delta^2\rho^2)}{\delta^2\rho^4}\int_0^{\infty}e^{-\delta\rho^2t}\cos\big(t\rho\sqrt{4-\delta^2\rho^2}\big)dt\\
\nonumber
&=\frac1{\delta\rho^2}-\frac{\rho^2(4-\delta^2\rho^2)}{\delta^2\rho^4}A,
\end{align}
which implies
$
A=\frac{\delta}{4}.
$
Substituting this $A$ into \eqref{int} yields
\begin{equation}\label{t}
\int_0^{\infty} e^{-\delta\rho^2t} \sin^2\big({\frac{t\rho\sqrt{4-\delta^2\rho^2}}{2}}\big)dt=\frac{4-\delta^2\rho^2}{8\delta\rho^2}.
\end{equation}
Then \eqref{f1} is bounded by
\begin{equation*}
\frac1{2\delta} \int_0^{2/\delta} {\rho^{2\alpha+2l+n-6}}d\rho<\infty
\end{equation*}
when 
\begin{equation}\label{con1}
2\alpha+2l+n-5>0 \quad \textrm{for all}\,\,\, l=0,1,2,\cdots.
\end{equation} 

To bound the last two integrals in \eqref{in}, we first calculate the integral in $t$:
\begin{equation*}
\int_0^{\infty}|m(t,\rho)|^2dt=\int_0^{\infty}\frac{1}{(\mu_1-\mu_2)^2}(e^{\mu_1 t}-e^{\mu_2 t})^2dt 
\end{equation*}
where $\mu_1$ and $\mu_2$ are the distinct real roots of the quadratic equation $\mu^2+\delta\rho^2\mu+\rho^2=0$ (see \eqref{mull}).
Since $\mu_1,\mu_2<0$, $\mu_1+\mu_2=-\delta\rho^2$ and $\mu_1\mu_2=\rho^2$, this integral is calculated as
\begin{align}\label{re}
\nonumber
\int_0^{\infty}|m(t,\rho)|^2dt&=\frac{1}{(\mu_1-\mu_2)^2}\int_0^{\infty}e^{2\mu_1t}-2e^{(\mu_1+\mu_2)t}+e^{2\mu_2t}dt \\
\nonumber
&=\frac{1}{(\mu_1-\mu_2)^2}\Big[\frac{1}{2\mu_1}e^{2\mu_1t} - \frac{2}{\mu_1+\mu_2}e^{(\mu_1+\mu_2)t}+\frac{1}{2\mu_2}e^{2\mu_2t}\Big]_0^{\infty}\\
&=-\frac{1}{2\mu_1\mu_2(\mu_1+\mu_2)}=\frac{1}{2\delta\rho^4}.
\end{align}
By using this and \eqref{Be}, the last two integrals in \eqref{in} are bounded by
\begin{align*}
\int_{2/\delta}^1 \rho^{2\alpha+2l+n-6}d\rho+\int_1^{\infty}\rho^{2\alpha-5}d\rho <\infty
\end{align*}
when 
\begin{equation}\label{con2}
\alpha<2.
\end{equation}

Finally, combining \eqref{con1} and \eqref{con2} yields the assumption $-\frac{n-5}{2}<\alpha<2$ in Lemma \ref{ker}.
(The condition $n\ge2$ follows immediately from this assumption.)

\subsection{The case $0<\delta\leq2$}
Since $2/\delta\ge1$ in this case, we split this time the integral with respect to $\rho$ in \eqref{K} into the regions $0<\rho<1$, $1<\delta<2/\delta$ and $\rho>2/\delta$:
\begin{align}\label{in2}
\nonumber
&\int_0^{1}\int_0^{\infty}|m(t,\rho)|^2\rho^{2\alpha}|J_{l+\frac{n-2}{2}}(\rho)|^2dtd\rho+\int_1^{2/\delta}\int_{0}^{\infty}|m(t,\rho)|^2\rho^{2\alpha}|J_{l+\frac{n-2}{2}}(\rho)|^2dtd\rho\\
&\qquad\qquad +\int_{2/\delta}^{\infty}\int_0^{\infty}|m(t,\rho)|^2\rho^{2\alpha}|J_{l+\frac{n-2}{2}}(\rho)|^2dtd\rho.
\end{align}

By using \eqref{Be} and \eqref{t}, the first two integrals in \eqref{in2} are bounded by
\begin{equation*}
\int_0^1 \rho^{2\alpha+2l+n-6}d\rho+\int_1^{\delta/2}  \rho^{2\alpha-5}d\rho <\infty
\end{equation*}
when $\alpha >-\frac{n-5}{2}$.
On the other hand, by using \eqref{Be} and \eqref{re}, the last integral in \eqref{in2} is bounded by 
\begin{equation*}
\int_{2/\delta}^{\infty} {\rho^{2\alpha-5}}d\rho <\infty
\end{equation*}
when $\alpha<2$.
\subsection*{Acknowledgment}
The authors would like to thank the anonymous referees for their valuable comments and helpful suggestions. 

\subsection*{Data availability}
No new data were created or analysed in this study.


\begin{thebibliography}{99}
	
\bibitem{CS}
R. W. Carroll and R. E. Showalter, \textit{Singular and Degenerate Cauchy Problems}, Mathematics in Science and Engineering , Vol. 127.
Academic Press [Harcourt Brace Jovanovich, Publishers], New York-London, 1976.

\bibitem{C}
J. B. Conway, \textit{A course in functional analysis}, Grad. Texts in Math., 96
Springer-Verlag, New York, 1985. xiv+404 pp.

\bibitem{ESS}
P. Elbau, O. Scherzer, and C. Shi, \textit{Singular values of the attenuated photoacoustic imaging operator}, J. Differential Equations 263 5330-5376 (2017).


\bibitem{WF}
D. Fang and C. Wang, \textit{Weighted Strichartz estimates with angular regularity and their applications}, Forum Math. 23 (2011), 181–205.

\bibitem{FR}
D. Finch and Rakesh, \textit{ The spherical mean value operator with centers on a sphere}, Inverse Problems, 23 (2007), S37-S49.

\bibitem{FPR}
D. Finch, S. Patch and Rakesh, \textit{Determining a function from its mean values over a family of spheres}, SIAM J. Math. Anal., 35 (2004), 1213-1240. 

\bibitem{Ga}
L. Grafakos, \textit{Classical Fourier Analysis}, 2nd edition, Graduate Texts in Mathematics, 249. Springer, New York, 2008.


\bibitem{G}
J. M. Greeberg, \textit{On the existence, uniqueness, and stability of solutions of the equation $\rho_0 X_tt=E(X_x)X_{xx}+\kappa X_{xxt}$}, J. Math. Anal. Appl. 25 (1969), 575-591.


\bibitem{kuchmentk08}
P.~Kuchment and L.~Kunyansky.
\newblock {\em Mathematics of thermoacoustic tomography.}
\newblock  European Journal of Applied Mathematics, 19:191--224, 2008.


\bibitem{kuchment14book}
P.~Kuchment.
\newblock {\em The Radon Transform and Medical Imaging}.
\newblock CBMS-NSF Regional Conference Series in Applied Mathematics. Society
  for Industrial and Applied Mathematics, 2014.

\bibitem{kowars10}
R.~Kowar and O.~Scherzer.
\newblock {\em Attenuation Models in Photoacoustics}.
\newblock Mathematical Modeling in Biomedical Imaging II, 2011. 


\bibitem{moonjmaa18}
S. Moon, \newblock {\em Inversion formulas and stability estimates of the wave operator on the hyperplane}. 
\newblock Journal of Mathematical Analysis and Applications, 466(2018), 490-497, 2018.
	
\bibitem{S1}	
R. T. Seeley, \textit{Spherical harmonics}, Amer. Math. Monthly 73 (1966), 115–121.


\bibitem{S}
R. E. Showalter, \textit{Regularization and approximation of second order evolution equations}, SIAM J. Math. Anal. 7 (1976), 461-472.

\bibitem{SW}
E. M. Stein and G. Weiss, \textit{Introduction to Fourier Analysis on Euclidean Spaces},
Princeton Mathematical Series 32, Princeton University Press, Princeton, N. J., 1971

\bibitem{W}
L. W. White, \textit{Approximation of point controls of second-order evolution equations of Sobolev type}, J. Math. Anal. Appl. 90 (1982), 117–126.

\bibitem{xuw05}
M. Xu, LV Wang.
\newblock {\em Universal back-projection algorithm for photoacoustic computed tomography.}
\newblock Physical Review E, 71 (1), 016706, 2005.

\bibitem{zangerl19}
		G. Zangerl, S. Moon, and M. Haltmeier, 
		\newblock {\em Photoacoustic tomography with direction dependent data: an exact series reconstruction approach.} 
		\newblock Inverse Problems, 35(11), 114005, (2019).

\end{thebibliography}
\end{document}